\documentclass{amsproc}
\usepackage{amsfonts}
\usepackage{amssymb}
\usepackage{amscd}
\usepackage{verbatim}
\usepackage{graphicx}
\usepackage{subfigure}
\usepackage{fancyhdr}
\usepackage[english]{babel}
\usepackage[T1]{fontenc}
\usepackage[latin1]{inputenc}
\usepackage{bm}

\input xy
\xyoption{all}

\newcommand{\mor}[3]{$\xymatrix@1@C=15pt{#3: #1\ar[r]& #2}$}
\newcommand{\appl}[2]{$\xymatrix@1@C=15pt{#1 \ar@{|->}[r]& #2}$}
\newcommand{\pf}{\medskip \noindent {\bf Proof :\ \ }}

\newcommand{\ol}[1]{\overline{#1}}

\newcommand{\cqfd}{\begin{flushright}$\Box$\end{flushright}}
\newcommand{\rad}[1]{{\rm rad}#1}

\begin{document}
\bibliographystyle{plain}

\begin{abstract}
We show that the graded commutative ring structure of the
Hochs\-child cohomology ${\rm HH}^*(A)$ is trivial in case $A$ is
a triangular quadratic string algebra. Moreover, in case $A$ is
gentle, the Lie algebra structure on ${\rm HH}^*(A)$ is also
trivial.
\end{abstract}

\title{The cohomology structure of string algebras}
\author[J.~C.~Bustamante ]{Juan Carlos Bustamante \footnote{
Hochschild cohomology rings, string algebras. }}
\address{Departamento de Matem\'atica - I.M.E, Universidade de S\~ao Paulo,
Caixa Postal 66281, 05315-970, S\~ao Paulo, Brazil}
\email{bustaman@ime.usp.br} \maketitle

\thispagestyle{empty}

\section*{Introduction}

Let $k$ be a field, and $A$ be a finite dimensional $k-$algebra.
In this situation, the Hochschild cohomology groups ${\rm
HH}^i(A,M)$ with coefficients in some $A-A-$bimodule $M$ can be
identified to the groups ${\rm Ext}_{A-A}^i(A, M)$. In case $_AM_A
=_AA_A$, we simply write ${\rm HH}^i(A)$. The sum ${\rm HH}^*(A) =
\coprod_{i\geq 0}{\rm HH}^i(A)$ is a graded commutative ring for
the Yoneda product, which coincides with a cup-product $\cup$.
Beside this, there is another product, namely the bracket product
$[-,-]$ which makes of ${\rm HH}^*(A)$ a graded Lie algebra.
Moreover, these two structures are related so that ${\rm HH}^*(A)$
is in fact a Gerstenhaber algebra \cite{G63}.

In general, very little information is known about these
structures. As Green and Solberg point out in \cite{GS03}, {\em
"the ring structure of ${\rm HH}^*(A)$ has often been observed to
be trivial"}. Although {\em "one knows that for many self
injective rings there are non zero products in ${\rm HH}^*(A)$"},
there are rather few known examples of algebras having finite
global dimension such that the ring structure of ${\rm HH}^*(A)$
is not trivial. This leads us to state the following:

\subsection*{Conjecture.} {\em Let $A=kQ/I$ be a monomial
triangular algebra. Then the ring structure of ${\rm HH}^*(A)$ is
trivial.}

\medskip

In case $A=kQ/I$ is a monomial algebra, one has a precise
description of a minimal resolution of  $_AA_A$. The
$A-A-$projective bimodules appearing in this resolution are
described in terms of paths of $Q$ (see \cite{B97}). In the
general case, this description leads to hard combinatorial
computations. However, in case $A$ is a monomial quadratic
algebra, the minimal resolution is particularly easy to handle.
Among monomial algebras, the class of string algebras, which are
always tame, is particularly well understood, at least from the
representation theoretic point of view \cite{BR87}. We study the
cohomology structure of triangular quadratic string algebras.

After recalling some notions concerning the products in ${\rm
HH}^*(A)$ in section \ref{sec:Prel}, we establish some technical
previous results in section \ref{sec:lemmata}. Finally, section
\ref{sec:results} is devoted to estate and show the main result of
this work, which proves the conjecture below in the particular
case of triangular quadratic string algebras.

\subsection*{Theorem  \ref{subsec:thm1} } {\em Let $A=kQ/I$ be a
triangular quadratic string algebra, $n>0$ and $m>0$. Then ${\rm
HH}^n(A) \cup {\rm HH}^m(A) =0 $.}

\medskip

In addition, we obtain a similar result concerning the Lie algebra
structure for a particular case of string algebras, the so-called
gentle algebras.

\subsection*{Theorem \ref{subsec:thm2}} {\em Let $A=kQ/I$ be a
triangular gentle algebra, $n>1$ and $m>1$. Then  $[{\rm
HH}^n(A),\ {\rm HH}^m(A)] =0$.}

\section{Preliminaries}\label{sec:Prel}

\subsection{Algebras} Let $Q=(Q_0, Q_1, s,t)$ be a finite quiver
(see \cite{BG82}), and $k$ be a commutative field. We consider
algebras of the form $kQ/I$, where $I$ an admissible ideal of the
path algebra $kQ$. This includes, for instance, all basic
connected and finite dimensional algebras over algebraically
closed fields (see \cite{BG82}).

While we briefly recall some particular concepts concerning bound
quivers and algebras, we refer the reader to \cite{BG82}, for
instance, for unexplained notions. The composition of two arrows
\mor{s(\alpha_1)}{t(\alpha_1)}{\alpha_1}, and
\mor{s(\alpha_2)}{t(\alpha_2)}{\alpha_2} such that $s(\alpha_2) =
t(\alpha_1)$ is the path $\alpha_1 \alpha_2$, and will be denoted
by $\xymatrix@C=35pt@1{\alpha_1 \alpha_2: s(\alpha_1)\ar@{~>}[r]&
t(\alpha_2)}$. If $Q$ has no oriented cycles, then $A$ is said to
be a {\em triangular algebra}. A two sided ideal $I\unlhd kQ$
generated by paths of $Q$ is said to be {\em monomial}. Moreover
if the generators of an admissible ideal $I$ are linear
combinations paths of length 2, then $I$ is said to be  {\em
quadratic}. An algebra $A=kQ/I$ is said to be {\em monomial}, or
{\em quadratic}, depending on whether $I$ is monomial or
quadratic. In the remaining part of this paper, all algebras are
triangular, quadratic and monomial.

Among monomial algebras, the so-called {\em string} algebras
\cite{BR87} are particularly well understood, at least from the
representation theoretic point of view. Recall from \cite{BR87}
that an algebra $A=kQ/I$ is said to be a {\em string algebra} if
$(Q,I)$ satisfies the following conditions :
\begin{enumerate}

\item [(S1)] $I$ is a monomial ideal,
\item [(S2)] Each vertex of $Q$ is the source and the target of at
most two arrows, and
\item [(S3)] For an arrow $\alpha$ in $Q$ there exists at most one
arrow $\beta$ and at most one arrow $\gamma$ such that $\alpha
\beta \not\in I$ and $\gamma \alpha\not\in I$.
\end{enumerate}

A string algebra $A=kQ/I$ is called a {\em gentle algebra}
\cite{AS87} if in addition $(Q,I)$ satisfies :

\begin{enumerate}
\item [(G1)] For an arrow $\alpha$ in $Q$ there exists at most one
arrow $\beta$ and at most one arrow $\gamma$ such that $\alpha
\beta \in I$ and $\gamma \alpha\in I$,
 \item [(G2)] $I$ is quadratic.
\end{enumerate}

\subsection{Hochschild cohomology} \label{subsec:Hochschild}
Given an algebra $A$ over a field $k$, the Hochs\-child cohomology
groups of $A$ with coefficients in some $A-A-$bimodule $M$,
denoted by ${\rm HH}^i(A,M)$,  are the groups ${\rm
Ext}^i_{A^e}(A,M)$ where $A^e = A \otimes_k A^{op}$ is the
enveloping algebra of $A$. In case $M$ is the $A-A-$bimodule ${}_A
A_A$, we simply denote them by ${\rm HH}^i(A)$. We refer the
reader to \cite{Hap89}, for instance, for general results about
Hochschild (co)-homology of algebras.

The lower cohomology groups, that is those for $0\leq i \leq 2$,
have clear interpretations and give information about, for
instance, the simple-connectedness or the rigidity properties of
$A$ (see \cite{AdlP96, G64}, for instance).

Beside this, the sum ${\rm HH}^*(A) = \coprod_{i\geq 0} {\rm
HH}^i(A)$ has an additional structure given by two products,
namely the cup-product $\cup$, and the bracket $[-,-]$, which are
defined in \cite{G63} at the cochain level. The definitions are
given in \ref{subsec:products}.

\subsection{Resolutions}\label{subsec:resolutions}

From \cite{C3}, we have a convenient projective resolution of $A$
over $A^e$ which is smaller that the standard Bar resolution. Let
$E$ be the (semi-simple) subalgebra of $A$ generated by $Q_0$. In
the remaining part of this note, tensor products will be taken
over $E$, unless it is explicitly otherwise stated. Note that, as
$E-E-$bimodules, we have $A\simeq E \oplus {\rm rad}A$. Let
$\rad{A}^{\otimes n}$ denote the $n^{th}$ tensor power of
$\rad{A}$ with itself. With these notations, one has a projective
resolution of $A$ as $A-A-$bimodule, which we denote by $\mathcal
{K}_{rad}^{\bullet}(A)$ $$\xymatrix@R=5pt{ \ldots \ar[r]& A
\otimes \rad{A}^{\otimes n}\otimes A \ar[r]^{b_{n-1}}& A \otimes
\rad{A}^{\otimes n-1}\otimes A \ar[r]^{\ \ \ \ \ \ \ \ \ \ \ \ \ \
b_{n-2}}&\ldots
\\
       \ar[r]& A \otimes \rad{A} \otimes A \ar[r]^{b_0}             &A\otimes A \ar[r]^{\epsilon}       & A\ar[r] &0 }$$
\noindent where $\epsilon$ is the multiplication of $A$, and
\begin{eqnarray*}
b_{n-1}(1\otimes r_1 \otimes \cdots \otimes r_n \otimes 1) & = &
r_1 \otimes r_2 \otimes \cdots \otimes r_n \otimes 1\\ &+&
\sum_{j=1}^{n-1}(-1)^j 1\otimes r_1 \otimes \cdots \otimes r_j r
_{j+1} \otimes \cdots \otimes r_n \otimes 1\\ &+& (-1)^n 1\otimes
r_1 \otimes \cdots \otimes r_n .
\end{eqnarray*}

\subsection*{Remark} Note that since the ideal $I$ is assumed to be monomial, ${\rm
rad}A$ is generated, as $E-E-$bimodule, by classes $\overline{p}=
p + I$ of paths of $Q$. Moreover, since the tensor products are
taken over $E$, then $A\otimes \rad{A}^{\otimes n} \otimes A$ is
generated, as $A-A-$bimodule, by elements of the form $1\otimes
\overline{p}_1 \otimes \cdots \otimes \overline{p}_n \otimes 1$
where $p_i$ are paths of $Q$ such that the ending point of $p_i$
is the starting point of $p_{i+1}$, for each $i$ such that $1\leq
i <n$.

\medskip
Keeping in mind that $I$ is quadratic, the minimal resolution of
Bardzell \cite{B97} has the following description: Let $\Gamma_0=
Q_0,\ \Gamma_1 = Q_1$, and for $n\geq 2$ let $\Gamma_n =
\{\alpha_1 \alpha_2\cdots \alpha_n |\ \alpha_i \alpha_{i+1} \in I,
\mbox{ for } 1\leq i < n\}$. For $n\geq 0$, we let $k\Gamma_n$ be
the $E-E-$bimodule generated by $\Gamma_n$. With these notations,
we have a minimal projective resolution of $A$ as $A-A-$module,
which we denote by $\mathcal {K}_{min}^{\bullet}(A)$:
$$\xymatrix@R=5pt{ \ldots \ar[r]& A \otimes k\Gamma_n \otimes A
\ar[r]^{\delta_{n-1}}& A \otimes k \Gamma_{n-1}\otimes A \ar[r]^{\
\ \ \ \ \ \ \ \ \  \delta_{n-2}}&\ldots
\\
\ar[r]& A \otimes k\Gamma_1 \otimes A \ar[r]^{\delta_0}&A\otimes
k\Gamma_0\otimes A \ar[r]^{\ \ \ \ \ \ \ \epsilon} & A\ar[r]&0 }$$

\noindent where, again, $\epsilon$ is the composition of the
isomorphism $A\otimes k\Gamma_0 \otimes A \simeq A \otimes A$ with
the multiplication of $A$, and, given $1\otimes \alpha_1 \cdots
\alpha_n\otimes 1 \in A\otimes k\Gamma_n \otimes A$ we have
$$\delta_{n-1}(1\otimes \alpha_1 \cdots \alpha_n\otimes 1) =
\alpha_1 \otimes \alpha_2 \cdots \alpha_n \otimes 1+ (-1)^n 1
\otimes \alpha_1 \cdots \alpha_{n-1} \otimes \alpha_{n}.$$

In order to compute the Hochschild cohomology groups of $A$, we
apply the functor ${\rm Hom}_{A^e}(-,A)$ to the resolutions
$\mathcal {K}_{min}^{\bullet}(A)$ and $\mathcal
{K}_{rad}^{\bullet}(A)$. To avoid cumbersome notations we will
write $\mathcal{C}_{rad}^n(A,A)$ instead of ${\rm
Hom}_{A^e}(A\otimes {\rm rad}A^{\otimes n} \otimes A, A)$ which in
addition we identify to ${\rm Hom}_{E^e}( {\rm rad}A^{\otimes n},
A)$. In a similar way we define $\mathcal{C}_{min}^n(A,A)$, and,
moreover, the obtained differentials will be denoted by
$b^\bullet$, and $\delta^\bullet$.

The products in the cohomology ${\rm HH}^*(A)$ are induced by
products defined using the Bar resolution of $A$ (see \cite{G63}).
Keeping in mind that our algebras are assumed to be triangular (so
that the elements of ${\mathcal C}^n_{min}(A,A)$ take values in
${\rm rad}A$), these products are easily carried to $\mathcal
{K}_{rad}^{\bullet}(A)$ (see section \ref{subsec:products}).
However, the spaces involved in that resolution are still "too
big" to work with. We wish to carry these products from $\mathcal
{K}_{rad}^{\bullet}(A)$ to $\mathcal {K}_{min}^{\bullet}(A)$. In
order to do so, we need {\em explicit} morphisms \mor{\mathcal
{K}_{min}^{\bullet}(A)}{\mathcal
{K}_{rad}^{\bullet}(A)}{\mu_\bullet}, and \mor{\mathcal
{K}_{rad}^{\bullet}(A)}{\mathcal
{K}_{min}^{\bullet}(A)}{\omega_\bullet} (compare with
\cite{Stra02} p.736).

\medskip
Define a morphism of $A-A-$bimodules \mor{A\otimes k\Gamma_n
\otimes A}{A \otimes {\rm rad}A^{\otimes n}\otimes A}{\mu_n} by
the rule $\mu_n(1\otimes \alpha_1 \cdots \alpha_n \otimes 1) = 1
\otimes \overline{\alpha}_1 \otimes \cdots \otimes
\overline{\alpha}_n \otimes 1.$ A straightforward computation
shows that $\mu_{\bullet} = (\mu_n)_{n\geq 0}$ is a morphism of
complexes, and that each $\mu_n$ splits.

On the other hand, an element  $1\otimes \overline{p}_1 \otimes
\cdots \otimes \overline{p}_n \otimes 1$ of $A\otimes
\rad{A}^{\otimes n} \otimes A$, can always be written as $1\otimes
\overline{p}'_1 \overline{\alpha}_1\otimes \overline{p}_2\otimes
\cdots\otimes \overline{\alpha}_n \overline{p}'_n\otimes 1$ with
$\alpha_1,\ \alpha_n$ arrows in $Q$. Define a map \mor{A \otimes
{\rm rad}A^{\otimes n}\otimes A}{A\otimes k\Gamma_n \otimes
A}{\omega_n} by the rule$$ \omega_n(1\otimes \overline{p}'_1
\overline{\alpha}_1\otimes
 \cdots\otimes \overline{\alpha}_n \overline{p}'_n\otimes 1)
= \left\{ \begin{array}{cl} \overline{p}'_1\otimes \alpha_1
p_2\cdots p_{n-1} \alpha_n \otimes \overline{p}'_n& \mbox{if }
\alpha_1 p_2\cdots p_{n-1} \alpha_n \in \Gamma_n\\ 0&
\mbox{otherwise.}
\end{array} \right.$$
Again, $\omega_{\bullet} = (\omega_n)_{n\geq 0}$ defines a
morphism of complexes, which is an inverse for $\mu_{\bullet}$.
Summarizing what precedes, we obtain the following lemma.

\subsection{Lemma}{\em With the above notations,
\begin{enumerate}
\item[a)] \mor{\mathcal {K}_{min}^{\bullet}(A)}{\mathcal
{K}_{rad}^{\bullet}(A)}{\mu_\bullet}, and \mor{\mathcal
{K}_{rad}^{\bullet}(A)}{\mathcal
{K}_{min}^{\bullet}(A)}{\omega_\bullet} are morphisms of
complexes, and
\item[b)] $\omega_\bullet \mu_\bullet = id$, so the morphisms are
quasi-isomorphisms.
\end{enumerate}}\cqfd

\subsection{Products in cohomology}\label{subsec:products} Following \cite{G63}, given
$f\in \mathcal{C}^n_{rad}(A,A),\ g\in \mathcal{C}^m_{rad}(A,A)$,
and $i\in \{1,\ldots, n\}$, define the element $f\circ_i g \in
\mathcal{C}^{n+m-1}_{rad}(A,A)$ by the rule $$f\circ_i g
(r_1\otimes \cdots \otimes r_{n+m-1}) = f(r_1\otimes \cdots
\otimes r_{i-1} \otimes g(r_i\otimes \cdots \otimes r_{i+m-1})
\otimes r_{i+m}\otimes \cdots \otimes r_{n+m-1}).$$ The
composition product $f\circ g$ is then defined as $\sum_{i=1}^n
(-1)^{(i-1)(m-1)} f\circ_i g.$ Let $f\circ g=0 $ in case $n=0$,
and define the {\em bracket} $$[f,g] = f\circ g -
(-1)^{(n-1)(m-1)} g\circ f.$$

On the other hand, the {\em cup-product} $f\cup g \in
\mathcal{C}^{n+m}_{rad}(A,A)$ is defined by the rule $$(f\cup
g)(r_1\otimes\cdots \otimes r_{n+m}) = f(r_1\otimes\cdots \otimes
r_n) g(r_{n+1}\otimes\cdots \otimes r_{n+m})$$

Recall that a Gerstenhaber algebra is a graded $k-$vector space
$A$ endowed with a product which makes $A$ into a graded
commutative algebra, and a bracket $[-,-]$ of degree $-1$ that
makes $A$ into a graded Lie algebra, and such that $[x,yz] =
[x,y]z + (-1)^{(|x|-1)|y|}y[x,z]$, that is, a graded analogous of
a Poisson algebra.

The cup product $\cup$ and the bracket $[-,-]$ define products
(still denoted $\cup$ and $[-,-]$) in the Hochschild cohomology
$\coprod_{i\geq 0} {\rm HH}^i(A)$, which becomes then a
Gerstenhaber algebra (see \cite{G63}).

Using $\mu_\bullet$ and $\omega_\bullet$ we can define analogous
products in $\mathcal{C}^\bullet_{min}(A,A)$. We study these
products later.

\section{Preparatory lemmata}\label{sec:lemmata} Recall that all our algebras are
assumed to be triangular, monomial and quadratic. The spaces
$\mathcal{C}^n_{min}(A,A) = {\rm Hom}_{E^e}(k\Gamma_n,A)$ have
natural bases ${\mathcal B}^n$ that we identify to the sets
$\{(\alpha_1\cdots \alpha_n, \overline{p})|\ \alpha_1\cdots
\alpha_n \in\Gamma_n,\ p$ is a path with $s(p)=s(\alpha_1),\
t(p)=t(\alpha_n),\ p\not\in I\}$. More precisely the map of
$E-E-$bimodules \mor{k\Gamma_n}{A}{f} corresponding to
$(\alpha_1\cdots \alpha_n, \overline{p})$ is defined by
$$f(\gamma_1\cdots \gamma_n) = \left\{\begin{array}{cl} \ol{p} &
\mbox{ if } \gamma_1\cdots \gamma_n = \alpha_1 \cdots \alpha_n,\\
0 & \mbox{ otherwise. } \end{array} \right.$$

We distinguish three different kinds of basis elements in
$\mathcal{C}^n_{min}(A,A)$, which yield three subspaces
$\mathcal{C}^n_{-},\ \mathcal{C}^n_{+}$, and $\mathcal{C}^n_0$ of
$\mathcal{C}^n_{min}(A,A)$:
\begin{enumerate}
\item $\mathcal{C}^n_{-}$ is generated by elements of the form $(\alpha_1\cdots \alpha_n,
\ol{\alpha}_1 \ol{p})$,
\item $\mathcal{C}^n_{0}$ is generated by elements of the form $(\alpha_1\cdots
\alpha_n,\ol{w})$ such that $w\not\in <\alpha_1 >$,  $w \not \in
<\alpha_n> $, and
\item $\mathcal{C}^n_{-}$ is generated by elements of the form $(\alpha_1\cdots
\alpha_n, \ol{q}\ol{\alpha}_n)$ such that $q\not\in <\alpha_1>$.
\end{enumerate}

Clearly, as vector spaces we have $\mathcal{C}^n_{min}(A,A) =
\mathcal{C}^n_{-} \amalg \mathcal{C}^n_{0} \amalg
\mathcal{C}^n_+$. Moreover, let $\mathcal{B}_{-},\
\mathcal{B}_{0}$, and $\mathcal{B}_{+}$ be the natural bases of
$\mathcal{C}^n_{-},\ \mathcal{C}^n_{+}$, and $\mathcal{C}^n_0$.

Consider an element $f=(\alpha_1\cdots \alpha_n, \ol{\alpha_1}
\ol{p}),\ \in \mathcal{B}^n_-$. First of all, since $\alpha_1
\alpha_2 \in I$, then $p\not \in <\alpha_2>$, and, since $A$
assumed to be triangular, then $p \not\in <\alpha_1>$. The
following figure illustrates this situation:
$$\xymatrix@C=25pt{e_0 \ar[r]^{\alpha_1} & e_1 \ar[r]^{\alpha_2}
\ar@{~>}`d[r]`[rrrr]_p[rrrr]& & \cdots & e_{n-1}
\ar[r]^{\alpha_n}& e_n\ar[r]^{\alpha_{n+1}}&e_{n+1}\\&&&&&&}$$

Now, define \mor{k\Gamma_n}{A}{f_+} in the following way:

$$f_+(\gamma_1\cdots \gamma_n) =\left\{\begin{array}{ll}
(-1)^{n+1} \ol{p\alpha_{n+1}} & \mbox{ if } \gamma_1\cdots
\gamma_n = \alpha_2 \cdots \alpha_{n+1} \in \Gamma_{n},\\ 0
&\mbox{ otherwise. } \end{array} \right.$$

Analogously, given $g=(\alpha_1\cdots \alpha_n,\ol{q \alpha_n})
\in \mathcal{B}^n_+$, we define \mor{k\Gamma_n}{A}{g_-} by

$$g_-(\gamma_1\cdots \gamma_n) =\left\{\begin{array}{ll} \ol{
\alpha_0 q} & \mbox{ if } \gamma_1\cdots \gamma_n = \alpha_0
\cdots \alpha_{n-1} \in \Gamma_{n},\\ 0 &\mbox{ otherwise. }
\end{array} \right.$$

It is easily seen that $f_+ \in \mathcal{C}^n_+$, and $g_-\in
\mathcal{C}^n_-$.

\subsection{Lemma}\label{subsec:lemma1}{\em Let $A=kQ/I$ be a
string triangular algebra,  $ f\in \mathcal{B}^n_-$, and $g \in
\mathcal{B}^n_+$.

\begin{enumerate}
\item [$a)$] $f-f_+ \in {\rm Im}\delta^{n-1}$, and
\item [$b)$] $g-g_- \in {\rm Im}\delta^{n-1}.$
\end{enumerate}
}

 \pf We only prove statement $a)$. Let $h=(\alpha_2\cdots
\alpha_n,\overline{p})\in \mathcal{C}_{min}^{n-1}(A,A)$. We show
that $f_+=f-\delta^{n-1} h$, indeed:

\begin{eqnarray*}
\delta^{n-1} h (\alpha_1\cdots \alpha_n) & = & \ol{\alpha}_1
h(\alpha_2\cdots \alpha_n) + (-1)^n h(\alpha_1\cdots
\alpha_{n-1})\ol{\alpha}_n\\
 & =& \ol{\alpha_1}\ol {p} \\& = &f(\alpha_1\cdots \alpha_n)
\end{eqnarray*}
\noindent thus $(f-\delta^{n-1} h)(\alpha_1\cdots \alpha_n) = 0 =
f_+(\alpha_1\cdots \alpha_n)$.

Now let $p=\beta p'$ with $\beta$ an arrow. Since $\alpha_1 p
\notin I$, then $\alpha_1 \beta \not\in I$. But $A$ is a string
algebra, thus condition (S3) ensures that for every arrow
$\alpha'_1$ such that $t(\alpha_1)=t(\alpha'_1)$ we have
$\alpha'_1\beta \in I$, thus $\alpha'_1 p \in I$.

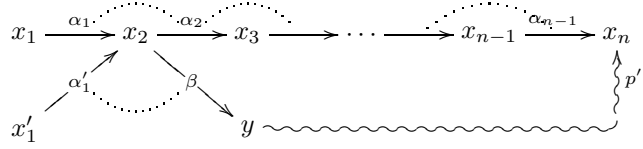
$$\xymatrix@C=10pt@R=7.5pt{x_1 \ar[rr]^{\alpha_1}&
\ar@/^1pc/@{.}[rr]&x_2\ar[rr]^{\alpha_2}\ar[ddrr]|{\beta}&\ar@/^1pc/@{.}[rr]
&x_3\ar[rr] & & \cdots \ar[rr] & \ar@/^1pc/@{.}[rr]&x_{n-1}
\ar[rr]^{\alpha_{n-1}}& & x_n \\&\ar@{.}@/_1pc/[rr]&&&&&&&&&\\
x'_1\ar[uurr]|{\alpha'_1}&&&&y\ar@{~>}`/3pt[uurrrrrr][uurrrrrr]_{p'}&&&&&&
}$$

Assume $\alpha'_1$ is such that $\alpha'_1\alpha_2\cdots
\alpha_n\in \Gamma_n$. In particular $\alpha'_1\alpha_2\in I$, so
that
\begin{eqnarray*}\delta^{n-1} h(\alpha'_1\alpha_2\cdots
\alpha_n)& = & \ol{\alpha'}_1 h(\alpha_2\cdots \alpha_n) + (-1)^n
h(\alpha'_1\cdots \alpha_{n-1})\ol{\alpha}_n\\ & = &
\ol{\alpha'}_1 \ol{p} \\ & = &0
\end{eqnarray*}
\noindent thus, again, $(f-\delta^{n-1} h)(\alpha'_1\cdots
\alpha_n) = 0 = f_+(\alpha'_1\cdots \alpha_n)$.

On the other hand, assume there exists an arrow $\alpha_{n+1}$
such that $\alpha_2\cdots \alpha_{n+1}\in \Gamma_n$. Then we have
\begin{eqnarray*}\delta^{n-1} h(\alpha_2\cdots \alpha_{n+1}) & = &
\ol{\alpha}_2 h(\alpha_3\cdots \alpha_{n+1}) +
(-1)^nh(\alpha_2\cdots \alpha_n)\\ & =& (-1)^n h(\alpha_2 \cdots
\alpha_n)\ol{\alpha}_{n+1}\\ &=& (-1)^n \ol{p \alpha_{n+1}}
\end{eqnarray*}
\noindent and hence
\begin{eqnarray*}
(f-\delta^{n-1} h)(\alpha_2 \cdots \alpha_{n+1})&= & 0 - (-1)^n
\ol{p \alpha_{n+1}}\\ & = & f_+(\alpha_2\cdots
\alpha_{n+1}).\end{eqnarray*} Finally, note that $f, f_+$, and $h$
vanish on any other element $\gamma_1\cdots \gamma_n$ of
$\Gamma_n$. \cqfd

According to the decomposition  $\mathcal{C}^n_{min}(A,A) =
\mathcal{B}^n_{-}\coprod \mathcal{B}^n_{0} \coprod
\mathcal{B}^n_{+}$, given an element  $\phi \in
\mathcal{C}^n_{min}(A,A)$ one can write $\phi = f + h + g$, where
$f\in\mathcal{B}^n_{-}, \ h \in  \mathcal{B}^n_{0}$, and $
g\in\mathcal{B}^n_{+}$. Define then $$ \phi_{\leq} = f_+ + h + g
\mbox{\ \ \ and\ \ \ }  \phi_{\geq}  =  f + h + g_-.$$

\subsection*{Remark} It follows from the preceding lemma that $\phi-\phi_{\leq}$, and
$\phi-\phi_{\geq}$ belong to ${\rm Im}\delta^{n-1}$. This will be
useful later.

\subsection{Lemma}\label{subsec:lemma2} {\em Let $A$ be a string triangular algebra, $\phi \in
\mathcal{C}^n_{min}(A,A)$, and $\alpha_1\cdots
\alpha_n\in\Gamma_n$. Then
\begin{enumerate}
\item [$a)$] $\phi_{\leq}(\alpha_1\cdots \alpha_n) = (W + Q \alpha_n) + I$ where
$W$ is a linear combination of paths none of which belongs to
$<\alpha_1>$ nor to $<\alpha_n>$, and $Q$ is a linear combination
of paths none of which belongs to $<\alpha_1>$, and

\item [$b)$] $\phi_{\geq}(\alpha_1\cdots \alpha_n) = (\alpha_1 P + W' ) + I$, where $W'$  is a
linear combination of paths none of which belongs to $<\alpha_1>$
nor to $<\alpha_n>$, and $P$ is a linear combination of paths none
of which belongs to $<\alpha_n>$.
\end{enumerate}}

\pf By construction we have $\phi_{\leq} \in \mathcal{B}^n_0
\coprod \mathcal{B}^n_+$, and $\phi_{\geq} \in \mathcal{B}^n_-
\coprod \mathcal{B}^n_0$. \cqfd

\subsection{Lemma}\label{subsec:lemma3}{\em Let $\phi\in{\rm Ker}
\delta^n$, and $\alpha_1\cdots \alpha_n \in \Gamma_n$ such that
$\phi(\alpha_1\cdots\alpha_n) \not=0$.
\begin{enumerate}
\item [$a)$] If there exists $\alpha_{n+1}$ such that $\alpha_1 \cdots
 \alpha_{n+1}\in \Gamma_{n+1}$ then $\phi_{\leq}(\alpha_1\cdots
\alpha_n)\ol{\alpha}_{n+1} = 0$,
\item [$b)$] If there exists $\alpha_0$ such that $\alpha_0 \cdots
\alpha_n\in \Gamma_{n+1}$ then $\ol{\alpha}_0 \phi_\geq
(\alpha_1\cdots \alpha_n) = 0$
\end{enumerate}}

\pf We only prove $a)$. Let $\phi\in{\rm Ker} \delta^n$. Then,
since $\phi-\phi_{\leq}\in {\rm Im}\delta^{n-1}$, we have
$\phi_{\leq} \in{\rm Ker} \delta^n$, thus
\begin{eqnarray*}
0 & = & \delta^{n}\phi_{\leq} (\alpha_1\cdots
\alpha_n\alpha_{n+1})\\
  & = & \ol{\alpha}_1 \phi_{\leq}(\alpha_2\cdots \alpha_{n+1}) +
  (-1)^{n+1} \phi_{\leq}(\alpha_1\cdots
  \alpha_n)\ol{\alpha}_{n+1}.
  \end{eqnarray*} The statement follows from the fact that
  paths are linearly independent, and from the preceding lemma.\cqfd

\section{The cohomology structure}\label{sec:results}

The morphisms \mor{\mathcal {K}_{min}^{\bullet}(A)}{\mathcal
{K}_{rad}^{\bullet}(A)}{\mu_\bullet}, and \mor{\mathcal
{K}_{rad}^{\bullet}(A)}{\mathcal
{K}_{min}^{\bullet}(A)}{\omega_\bullet} allow us to carry the
products defined in ${\mathcal C}_{rad}^{\bullet}(A,A)$ to
${\mathcal C}_{min}^\bullet(A,A)$. In this way we obtain a cup
product and a bracket, which we still denote $\cup$ and $[-,-]$.
More precisely, applying the functor ${\rm Hom}_{A^e}(-,A)$ and
making the identifications of section 1.5, we obtain morphisms of
complexes \mor{\mathcal {C}_{min}^{\bullet}(A,A)}{\mathcal
{C}_{rad}^{\bullet}(A,A)}{\mu^\bullet}, and \mor{\mathcal
{C}_{rad}^{\bullet}(A,A)}{\mathcal
{C}_{min}^{\bullet}(A,A)}{\omega^\bullet}. Given $f\in
\mathcal{C}^n_{min}(A,A),\ g\in \mathcal{C}^m_{min}(A,A)$, define
$f\cup g \in\mathcal{C}^{n+m}_{min}(A,A)$ as $$f\cup g =
\mu^{n+m}(\omega^n f \cup \omega^m f).$$ Thus, given an element
$\alpha_1\cdots \alpha_n\beta_1\cdots \beta_m \in \Gamma_{n+m}$,
we have $$f\cup g (\alpha_1\cdots \alpha_n \beta_1\cdots \beta_m)=
f(\alpha_1 \cdots \alpha_n)g(\beta_1\cdots \beta_m).$$ As usual,
we have $\delta^{n+m}(f\cup g)\ =\ \delta^n f\cup g\ +\ (-1)^m
f\cup \delta^m g$, so the product $\cup$ defined in ${\mathcal
C}_{min}^\bullet(A,A)$ induces a product at the cohomology level.
In fact, the latter coincides with the cup-product of section
\ref{subsec:products} and the Yoneda product. In what follows, we
work with the product  $\cup$ defined in ${\mathcal
C}_{min}^\bullet(A,A)$.

\subsection{Theorem}\label{subsec:thm1}{\em Let $A=kQ/I$ be a
triangular quadratic string algebra, $n>0$ and $m>0$. Then ${\rm
HH}^n(A) \cup {\rm HH}^m(A) =0 $.}

\pf Let $\ol{\phi}\in {\rm HH}^n(A)$, and $\ol{\psi}\in {\rm
HH}^m(A)$. We will show that $\ol{\phi_{\leq}}\cup \ol{\psi_\geq}
= 0$ in ${\rm HH}^{n+m}(A)$. Assume the contrary is true. In
particular there exist an element $\alpha_1\cdots \alpha_n\beta_1
\cdots \beta_m\in \Gamma_{n+m}$ such that
$$\phi_\leq(\alpha_1\cdots \alpha_n) \psi_\geq (\beta_1 \cdots
\beta_m) \not= 0.$$  With the notations of lemma
\ref{subsec:lemma2}, this is equivalent to
\begin{equation}\label{eqn:eq1}
(W+Q \alpha_n)(W'+ \beta_1 P)\not \in I. \end{equation} Now, since
$\alpha_n\beta_1\in I$, lemma \ref{subsec:lemma3} gives
$$(W+Q\alpha_n)\beta_1 \in I, \mbox{\ \  and\ \  }
\alpha_1(W'+\beta_1 P) \in I.$$

Thus, equation (\ref{eqn:eq1}) gives $WW'\not\in I$. Let $w$ and
$w'$ be paths appearing in $W$ and $W'$ respectively, such that
$ww'\not \in I$. Moreover, write $w=u\gamma,\ w' = \gamma'u'$ with
$\gamma, \gamma'$ arrows of $Q$. Since $ww'\not \in  I$, we have
that $\gamma \gamma'\not\in I$. But then, since $(Q,I)$ satisfies
(S3), this implies $\alpha_n \gamma'\in I$, so that $\alpha_1
\cdots \alpha_n\gamma' \in \Gamma_{n+1}$ and, again, lemma
\ref{subsec:lemma3} gives $(W +Q\alpha_n) \gamma' \in I$ so that
$W \gamma' \in I$, and, in particular $w\gamma' \in I$, a
contradiction.\ \ $\Box$

\subsection*{Remark} At this point, it is important to note that
even if  lemmas \ref{subsec:lemma1}, \ref{subsec:lemma2} and
\ref{subsec:lemma3} were stated assuming that $A=kQ/I$ is a string
triangular algebra, the condition $(S2)$ in the definition of
string algebras has never been used. Thus, the preceding theorem
holds for every monomial quadratic, triangular algebra $A=kQ/I$
such that $(Q,I)$ verifies $(S3)$. This includes, for instance,
the gentle algebras.

\medskip We now turn our attention to the Lie algebra structure ${\rm HH}^*(A)$.

\medskip Given $\phi\in \mathcal{C}^n_{min}(A,A),\ \psi\in
\mathcal{C}^m_{min}(A,A)$, and $i\in\{1,\ldots ,n\}$ we define
$\phi\circ_i \psi \in \mathcal{C}^{n+m-1}_{min}(A,A)$ as
$$\phi\circ_i \psi = \mu^{n+m-1}\left( (\omega^n \phi) \circ_i
(\omega^m \psi) \right)$$ and from this, the bracket $[-,-]$ is
defined as in \ref{subsec:products}. In particular, and this is
the crucial point for what follows, given $\alpha_1 \cdots
\alpha_{n+m-1} \in \Gamma_{n+m-1}$ we have $$\phi\circ_i \psi
(\alpha_1 \cdots \alpha_{n+m-1}) = \left\{\begin{array}{ll}
\phi(\alpha_1\cdots \alpha_{i-1}\psi(\alpha_i \cdots
\alpha_{i+m-1}) \alpha_{i+m}\cdots \alpha_{n+m-1})& \mbox{ if } \\
\ \ \ \ \ \alpha_1\cdots \alpha_{i-1}\psi(\alpha_i \cdots
\alpha_{i+m-1}) \alpha_{i+m}\cdots \alpha_{n+m-1}& \in \Gamma_n,\\
0 \mbox{\ \  otherwise.}&\end{array}\right.$$

\medskip Again, one can verify that $\delta^{n+m-1}[f,g] =
(-1)^{m-1}[\delta^n f, g]\ +\ [f,\delta^m g]$, so that $[-,-]$
induces a bracket, which we still denote by $[-,-]$ at the
cohomology level.

\medskip This leads us to the following result

\subsection{Theorem}\label{subsec:thm2} {\em Let $A=kQ/I$ be a triangular gentle
algebra. Then, for $n>1$ and $m>1$, we have $[{\rm HH}^n(A), {\rm
HH}^m(A)] = 0$.}

\pf We show that, in fact, under the hypothesis, the products
$\circ_i$ are equal to zero at the cochain level. This follows
immediately from the discussion above. Indeed, with that
notations, if $\alpha_1 \cdots \alpha_{n+m-1} \in \Gamma_{n+m-1}$
we have, in particular that $\alpha_{i-1}\alpha_i \in I$. But $A$
being gentle, there is no other arrow $\beta$ in $Q$  with
$s(\beta) = t(\alpha_{i-1})$ such that $\alpha_{i-1}\beta \in
I$.\cqfd

The following example shows that the previous result can not be
extended to string algebras which are not gentle.

\subsection{Example}Let $A=kQ/I$ where $Q$ is the quiver  $$
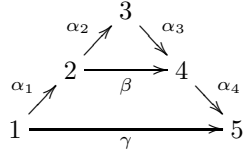
\xymatrix@C=10pt@R=10pt{    &                 &
3\ar[dr]^{\alpha_3}& &\\
                &2\ar[ur]^{\alpha_2} \ar[rr]_\beta &         & 4 \ar[dr]^{\alpha_4}   &\\
1\ar[ur]^{\alpha_1}\ar[rrrr]_\gamma&                &         &
&5}$$ and $I$ is the ideal generated by paths of length 2. This is
a string algebra which is not gentle. From theorem 3.1 in
\cite{C98} we get$${\rm dim}_k {\rm HH}^i(A) = \left\{
\begin{array}{ll} 1 & \mbox{ if } i\in\{0,2,3,4\}, \\
                  2 & \mbox{ if } i=1, \\
                  0 & \mbox{ otherwise. }\end{array} \right.$$
Keeping in mind the identifications of section 2, the generators
of ${\rm HH}^i(A)$ are the elements corresponding to $(\beta,
\bar{\beta})$, and $(\gamma, \bar{\gamma})$ for $i=1$;
$(\alpha_2\alpha_3,\bar{\beta})$ for $i=2$; $(\alpha_1 \beta
\alpha_4, \bar{\gamma})$, for $i=3$; and
$(\alpha_1\alpha_2\alpha_3\alpha_4, \bar{\gamma})$ for $i=4$. A
 straightforward computation shows that $[{\rm
                  HH}^n(A), {\rm HH}^m(A)] = {\rm HH}^{n+m-1}(A)$.

\medskip
The first Hochschild cohomology group of an algebra $A$ is by its
own right a Lie algebra. In case the algebra $A$ is monomial, this
structure has been studied in \cite{Stra02}. The following result
gives information about the role played by ${\rm HH}^1(A)$ in the
whole Lie algebra ${\rm HH}^*(A)$ in our context.

\subsection{Proposition}\label{subsec:propo}{\em Let $A=kQ/I$ be a
triangular monomial and quadratic algebra. Then $[{\rm HH}^n(A),
{\rm HH}^1(A)] = {\rm HH}^n(A)$, whenever $n>1$.}

\pf In fact, for every $f\in {\mathcal C}^n_{min}(A,A)$ there
exists $g\in {\mathcal C}^1_{min}(A,A)$ such that $[f,g] = f$, and
$\overline{g}\not=0$ in ${\rm HH}^1(A)$. Clearly it is enough to
consider an arbitrary basis element  $f\in {\mathcal
C}^n_{min}(A,A)$. Let $f$ be such an element, corresponding to
$(\alpha_1\cdots \alpha_n,\ol{p})$. Define $g\in{\mathcal
C}^1_{min}(A,A)$  by $$g(\gamma) =\left\{
\begin{array}{ll} \alpha_1 & \mbox{ if } \gamma = \alpha_1, \\
                  0& \mbox{ otherwise }\end{array} \right.$$
It is straightforward to verify that $f\circ g = f\circ_1 g$, and,
since we assume $A$ triangular and $n>1$, that  $g\circ f = 0$ so
that $[f,g]=f$. Moreover, direct computations show that $\bar{g}
\not=0$ in ${\rm HH}^1(A)$.\cqfd

\subsection*{Acknowledgements} We would like to thank professor E.
N. Marcos for several discussions and comments, as well as
professors E. N. Marcos and I. Assem for carefully reading a
preliminary version of this work and several useful suggestions.
The author gratefully acknowledges the financial support from
F.A.P.E.S.P., Brazil.

\bibliography{biblio}
\end{document}